\documentclass[12pt]{article}
\usepackage{amsmath}
 \usepackage{amssymb}
 \usepackage{array}
 \usepackage{geometry}
 \usepackage{graphicx}
 \usepackage{enumerate}
\newtheorem {thm}{Theorem}[section]
\newtheorem{lem}{Lemma}[section]

\newtheorem {ex}{Example}[section]
\newtheorem {de}{Definition}[section]
\numberwithin{equation}{section} 
\begin{document}
\begin{center}
{\Large \textbf{One-sided L\'{e}vy stable distributions 
}}%
\end{center}

{\begin{center} Jung Hun Han   \footnote{ Corresponding Author
Address: Centre for Mathematical Sciences, Pala Campus, Arunapuram
P.O., Pala, Kerala-686 574, India, ~~Email : jhan176@yahoo.com,
jhan176@gmail.com }
\end{center}}

%
\begin{abstract}
In this paper, we show new representations of one-sided L\'{e}vy
stable distributions for irrational L\'{e}vy indices of the type
$\left(\frac{p}{q}\right)^{\frac{l_{2}}{l_{1}}}$  which are not
covered in \cite{pg1} : for rational L\'{e}vy indices.
Furthermore, other equivalent representations for a distribution
of a rational L\'{e}vy index is described. We also give a simplest
proof for the formulae which cover the cases for rational L\'{e}vy
indices. Finally we introduce the concepts of L\'{e}vy smashing
and L\'{e}vy-smashed gamma stochastic processes.
\end{abstract}

It was written in August, 10, 2010.
%
%
%
%
%
\section{\bf  Introduction and preliminary results}\label{}
In \cite{pg1}, Penson and G\'{o}rska obtained a general form of
one-sided L\'{e}vy stable distributions expressed as a Meijer
$G$-function for rational L\'{e}vy indices by putting $v=0$ and
$a=1$ in the formula 2.2.1.19 in vol. 5 of \cite{pbm1}. 
In \cite{h1}, the role of Mathai transformation in the theory of
fractional calculus, which connects ordinary integral to
fractional integral through their kernels, is described and
$\alpha$-level space is defined. Furthermore it is insisted that
fractional integral and derivative preserve the L\'{e}vy structure
defined in \cite{h1}. The L\'{e}vy structure is nothing but the
integrand $\frac{\Gamma(\frac{1}{\alpha} - \frac{s}{\alpha}
)}{\alpha \Gamma(1-s) }$ of the $H$-function representation of the
L\'{e}vy distribution with $\alpha$ known as the L\'{e}vy index
which lies between $0$ and $1$. Hence the case of simple
irrational L\'{e}vy indices of the type $\left(\frac{p}{q}
\right)^{\frac{l_{2}}{l_{1}}}$ can be handled.\\

 In this paper, we provide formulae for one-sided L\'{e}vy distribution of
irrational L\'{e}vy index for $0<\alpha<1$ and some techniques to
generate infinitely many new representations of one-sided
L\'{e}vy distribution. Furthermore, we introduce the concept of
\textit{L\'{e}vy smashing} as a consequence of L\'{e}vy effect on the family of gamma density functions and
L\'{e}vy-smashed stochastic processes.\\

We will use the following integral representation of the gamma
function:
\begin{eqnarray}
  \nonumber\Gamma (z) &=& p^{z}\int^{\infty}_{0}t^{z-1}e^{-pt}dt, ~~\mathfrak{R}(p)>0, \mathfrak{R}(z)>0 \\
\nonumber   &=&\lim_{n\rightarrow\infty} \frac{n!~
n^{z}}{z(z+1)\cdots (z+n)},~~z\neq 0, 1, 2, 3, \cdots.
   \end{eqnarray}
Pochammer symbol is defined as
\begin{eqnarray}\label{} \nonumber(b)_{k}&=&b(b+1)\cdots (b+k-1), ~(b)_{0}=1,~ b \neq 0\\
\nonumber&=& (a)_{k}= \frac{\Gamma(a+k)}{\Gamma(a)} \mbox{
whenever the gammas exist.}
 \end{eqnarray}
%
For $H$-function representations and their convergence conditions,
consult with \cite{mh3}, \cite{mhs2}:
\begin{equation}
\nonumber H^{m,n}_{p,q}\left[z \big|^{\left(a_{1},A_{1}\right),
\left(a_{2},A_{2}\right),\cdots,\left(
a_{p},A_{p}\right)}_{(b_{1},B_{1}),\left((b_{2},B_{2}
\right),\cdots,\left((b_{q},B_{q}\right) }\right]
=\frac{1}{2 \pi i }\oint_{L} \frac{\{\prod^{m}_{j=1}
\Gamma(b_{j}+B_{j}s)\}\{ \prod^{n}_{j=1}
\Gamma\left(1-a_{j}-A_{j}s\right)\} }{ \{\prod^{q}_{j=m+1}
\Gamma\left(1-b_{j}-B_{j} s\right)\}\{
\prod^{p}_{j=n+1}\Gamma\left( a_{j}+A_{j} s\right)\} } z^{ -s} ds.
\end{equation}
Generalized Wright function, which is well explained in
\cite{mhs2} and which is a particular case of the $H$-function, is
\begin{equation}
\nonumber~_{p}\Psi_{q}\left[z \big|^{\left(a_{1},A_{1}\right),
\left(a_{2},A_{2}\right),\cdots,\left(
a_{p},A_{p}\right)}_{(b_{1},B_{1}),\left(b_{2},B_{2}
\right),\cdots,\left(b_{q},B_{q}\right) }\right]
=\sum_{n=0}^{\infty} \frac{\prod^{p}_{i=1} \Gamma\left(
a_{i}+A_{i} n\right)  }{ \prod^{q}_{j=1}\Gamma(b_{j}+B_{j}n)~n!
 } z^{n}
\end{equation}
where $a_{i},b_{j}\in \mathbb{C}$ and $A_{i},B_{j}\in \mathbb{R}$:
$A_{i}\neq 0,B_{j} \neq 0, i=1,\cdots, p, j=1,\cdots,q$;
$\sum_{j=1}^{q}B_{j}-\sum_{i=1}^{p}A_{i}>-1$ for absolute
convergence.\\
Hypergeometric series, which is well explained in \cite{mh3} and
which is a particular case of the $H$-function, is
\begin{equation}
\nonumber~_{p}F_{q}\left( a_{1},\cdots, a_{p};b_{1},\cdots,b_{q}
;z\right)
=\sum_{n=0}^{\infty} \frac{\prod^{p}_{i=1} \left( a_{i}\right) }{
\prod^{q}_{j=1}(b_{j})~n!
 } z^{n}
\end{equation}
which is absolutely convergent for all $z$ in $\mathbb{C}$ if
$p\leq q$.
\begin{de}
L\'{e}vy jump function is defined as follows
\begin{equation}
\nonumber\left[\frac{j}{q}\right]_{n}=\begin{cases} \frac{j}{q}, ~~~j<n,~j,n\in\mathbb{Z^{+}}\\
\frac{j+1}{q},~~~j\geq n,~j,n\in\mathbb{Z^{+}}. \end{cases}
\end{equation}
\end{de}
\begin{lem}\label{}
[ Stirling asymptotic formula ] \cite{ mh3}\\
For $|z|\rightarrow \infty$ and $\alpha$ a bounded quantity,
\begin{equation}\label{} \Gamma (z+ \alpha) \approx (2\pi)^{1/2}z^{z+\alpha
-1/2}e^{-z}.
\end{equation}
\end{lem}
\begin{lem}\label{h9}
\begin{eqnarray*}
\Gamma\left(-v+\frac{i}{q}-\left[\frac{j}{q}\right]_{i}\right)
&=&\frac{(-1)^{v}\Gamma\left(\frac{i}{q}-\left[\frac{j}{q}\right]_{i}\right)}
{\left(1-\frac{i}{q}+\left[\frac{j}{q}\right]_{i}\right)_{v}}\\
\Gamma\left(-v+\frac{iq-jp}{pq}\right)&=&\frac{(-1)^{v}\Gamma\left(\frac{iq-jp}{pq}
   \right)}
{\left(1-\frac{iq-jp}{pq}\right)_{v}}
\end{eqnarray*}
and
\begin{equation}
 \nonumber H^{m,n}_{p,q}\left[z
 \big|^{(a_{p},A_{p})}_{(b_{q},B_{q})}\right]=
 k~H^{m,n}_{p,q}\left[z^{k}
 \big|^{(a_{p},kA_{p})}_{(b_{q},kB_{q})}\right]
 \end{equation}
where $~v,n,p,q,i,j\in\mathbb{Z^{+}}$ and $k>0$.
 \end{lem}
Proof. Use the formula in \cite{mh3}, p44,
$\Gamma(\beta+1-v)=\frac{(-1)^{v}\Gamma(\beta+1)}{(-\beta)_{v}}$
and the properties in \cite{mhs2}, which complete the proof.
\section{\bf L\'{e}vy stable distribution of  L\'{e}vy index
$\left(\frac{p}{q}\right)^{\frac{l_{2}}{l_{1}}}$}\label{}

 Let
\begin{equation}\label{}
f_{\alpha}(x)=\frac{1}{2 \pi i}\oint_{L}
\frac{\Gamma(\frac{1}{\alpha} - \frac{s_{1}}{\alpha}  )}{\alpha
\Gamma(1-s_{1}) }x^{-s_{1}}ds_{1}, ~~ 1>\alpha>Re(s_{1})>0.
\end{equation}
$\frac{\Gamma(\frac{1}{\alpha} - \frac{s}{\alpha} )}{\alpha
\Gamma(1-s) }$ appears firstly in \cite{ma1} in the literature.
Then $f_{\alpha}(x)$ is the well known Levy density function
having the Laplace transform $e^{-t^{\alpha}}$. 

\begin{thm}
Let $0<\alpha<1$ and $p,q, l_{1},l_{2}$ be arbitrary positive
integers such that
$0<\alpha=\left(\frac{p}{q}\right)^{\frac{l_{2}}{l_{1}}}<1,~~p<q$
and $l=\left(\frac{p}{q} \right)^{\frac{l_{2}}{l_{1}}-1} $.\\

If $l$ is not a positive integer and $l_{1}\neq l_{2}$, then
\begin{eqnarray}
\nonumber&&f_{(\alpha,l)}(x)=\frac{l\sqrt{pq}}{x(2
\pi)^{\frac{q-p}{2}}}~H^{q,0}_{p,q}\left[\frac{p^{pl }}{x^{pl
}q^{q}}
 \big|^{\left(1,l\right),\left(\frac{1}{p},l\right),\cdots,\left(
\frac{p-1}{p},l\right)}_{(1,1),\left(\frac{1}{q},1
\right),\cdots,\left(\frac{q-1}{q},1\right) }\right]\\
%
\nonumber &&=\frac{l\sqrt{pq}}{x(2
\pi)^{\frac{q-p}{2}}}\left(\frac{p^{lp}}{x^{lp}
q^{q}}\right)~_{q-1}\Psi_{p}\left[-\frac{p^{pl}}{x^{pl}q^{q}}
{\Huge |}^{(-\frac{q-1}{q},-1)\cdots
 (-\frac{2}{q},-1) (-\frac{1}{q},-1)}_{(1-l,-l)
(\frac{q-plq}{pq},-l)\cdots (\frac{(p-1)q-plq}{pq},-l)}\right]
 +\frac{l\sqrt{pq}}{x(2
\pi)^{\frac{q-p}{2}}}\sum^{q-1}_{j=1}
\left(\frac{p^{jl\frac{p}{q}}}{x^{jl\frac{p}{q}} q^{j}}\right)\\
&&\times ~_{q-1}\Psi_{p}\left[-\frac{p^{pl}}{x^{pl}q^{q}} {\Huge
|}^{\left(1-\left[\frac{j}{q}\right]_{q},-1\right)\left(\frac{2}{q}-\left[\frac{j}{q}\right]_{2},-1\right)
\left(\frac{3}{q}-\left[\frac{j}{q}\right]_{3},-1\right)\cdots
\left(\frac{q-1}{q}-\left[\frac{j}{q}\right]_{q-1},-1\right)
 }_{\left(\frac{pq-j
pl}{pq},-l\right)\left(\frac{q-j pl}{pq},-l\right)\cdots
\left(\frac{(p-1)q-j pl}{pq},-l\right)}\right].\label{h1}
\end{eqnarray}

If $l=1$ and $l_{1}= l_{2}$, then
\begin{eqnarray}
\nonumber&& f_{(\alpha,1)}(x)=\frac{\sqrt{pq}}{x(2
\pi)^{\frac{q-p}{2}}}~H^{q-1,0}_{p-1,q-1}\left[\frac{p^{p }}{x^{p
}q^{q}}
 \big|^{\left(\frac{1}{p},1\right),\cdots,\left(
\frac{p-1}{p},1\right)}_{\left(\frac{1}{q},1
\right),\cdots,\left(\frac{q-1}{q},1\right) }\right]\\
%
%
\nonumber&&=\frac{\sqrt{pq}}{x(2
\pi)^{\frac{q-p}{2}}}\sum^{q-1}_{j=1}
\left(\frac{p^{j\frac{p}{q}}}{x^{j\frac{p}{q}} q^{j}}\right)
\frac{
\prod_{i_{1}=2}^{q-1}\Gamma\left(\frac{i_{1}}{q}-\left[\frac{j}{q}\right]_{i_{1}}\right)
} {\prod_{i_{2}=1}^{p-1}\Gamma\left(\frac{i_{2}q-j p}{pq}\right)}
\\
&&\times ~_{p-1}F_{q-2}\left[(-1)^{q-p}\frac{p^{p}}{x^{p}q^{q}}
{\Huge |}_{\left(1-\frac{2}{q}+\left[\frac{j}{q}\right]_{2}\right)
\left(1-\frac{3}{q}+\left[\frac{j}{q}\right]_{3}\right)\cdots
\left(1-\frac{q-1}{q}+\left[\frac{j}{q}\right]_{q-1}\right) }^{
\left(1-\frac{q-j p}{pq}\right)\cdots \left(1-\frac{(p-1)q-j
p}{pq}\right)}\right].\label{h2}
\end{eqnarray}

If $l$ belongs to the set $\{2,3,4,\cdots\}$ and $l_{1}\neq
l_{2}$, then
\begin{eqnarray}
\nonumber&&f_{(\alpha,l)}(x)=\frac{\sqrt{pql}}{x(2
\pi)^{\frac{q+1-p-l}{2}}}~H^{q-1,0}_{p+l-2,q-1}\left[\frac{p^{pl
}l^{l}}{x^{pl }q^{q}}
 \big|^{\left(\frac{1}{l},1\right),\cdots,\left(
\frac{l-1}{l},1\right),\left(\frac{1}{p},l\right),\cdots,\left(
\frac{p-1}{p},l\right)}_{\left(\frac{1}{q},1
\right),\cdots,\left(\frac{q-1}{q},1\right) }\right]\\
\nonumber&& =\frac{\sqrt{pql}}{x(2
\pi)^{\frac{q+1-l-p}{2}}}\sum^{q-1}_{j=1}
\left(\frac{p^{jl\frac{p}{q}}l^{j\frac{l}{q}}}{x^{jl\frac{p}{q}} q^{j}}\right)\\
&&\times
~_{q-2}\Psi_{p-1+l-1}\left[-\frac{p^{pl}l^{l}}{x^{pl}q^{q}} {\Huge
|}^{\left(\frac{2}{q}-\left[\frac{j}{q}\right]_{2},-1\right)
\left(\frac{3}{q}-\left[\frac{j}{q}\right]_{3},-1\right)\cdots
\left(\frac{q-1}{q}-\left[\frac{j}{q}\right]_{q-1},-1\right) }_{
\left(\frac{q-j l}{lq},-1\right)\cdots \left(\frac{(l-1)q-j
l}{lq},-1\right)\left(\frac{q-j pl}{pq},-l\right)\cdots
\left(\frac{(p-1)q-j pl}{pq},-l\right)}\right].\label{h3}
\end{eqnarray}

\end{thm}

Proof. Here we use a transformation
 $1-s_{1}=\alpha s$, then
\begin{equation}\label{h4}
 f_{\alpha}(x)=\frac{1}{2 \pi i}\oint_{L} \frac{\Gamma(s
)}{ \Gamma(\alpha s) }x^{\alpha s-1}ds, ~~
\frac{1}{\alpha}>Re(s)>\frac{1-\alpha}{\alpha}>0.
\end{equation}
Now, by using the Gauss-Legendre formula (Multiplicative formula)
for gamma function, namely
\begin{eqnarray*}
\nonumber\Gamma(mz)&=&(2\pi
)^{\frac{1-m}{2}}m^{mz-\frac{1}{2}}\Gamma(z)\Gamma(z+\frac{1}{m})\cdots\Gamma(z+\frac{m-1}{m}),~
m=1,2,\cdots,
\end{eqnarray*}
we have
\begin{eqnarray*}\label{}
\Gamma(\alpha s)&=& \Gamma\left(\left(\frac{p}{q}
\right)^{\frac{l_{2}}{l_{1}}} s\right)
=\Gamma\left(p\frac{ls}{q}\right)=(2 \pi)^{\frac{1-p}{2}}
p^{p\frac{ls}{q} -\frac{1}{2}}\Gamma\left(\frac{ls}{q}\right)
\Gamma\left(\frac{ls}{q}+\frac{1}{p}\right) \cdots
\Gamma\left(\frac{ls}{q}+\frac{p-1}{p}\right).
\end{eqnarray*}
Apply to the integrand in \eqref{h4}, then
\begin{eqnarray}\label{}
\nonumber f_{(\alpha,l)}(x)=\frac{1}{x}\frac{1}{2 \pi i }\oint_{L}
\frac{\Gamma(s )~x^{\frac{p}{q} ls}}{ (2 \pi)^{\frac{1-p}{2}}
p^{p\frac{ls}{q} -\frac{1}{2}}\Gamma(\frac{ls}{q})
\Gamma(\frac{ls}{q}+\frac{1}{p}) \cdots
\Gamma(\frac{ls}{q}+\frac{p-1}{p})}ds
\end{eqnarray}
where $\alpha=\frac{p}{q} l=\frac{p}{q}\left(\frac{p}{q}
\right)^{\frac{l_{2}}{l_{1}}-1}=\left(\frac{p}{q}
\right)^{\frac{l_{2}}{l_{1}}} .$
\\We use $$H^{m,n}_{p,q}\left[z
 \big|^{(a_{p},A_{p})}_{(b_{q},B_{q})}\right]=
 k~H^{m,n}_{p,q}\left[z^{k}
 \big|^{(a_{p},kA_{p})}_{(b_{q},kB_{q})}\right]$$ in Lemma \ref{h9}.
 Then
\begin{eqnarray*}\label{}
&&f_{(\alpha,l)}(x)=\frac{q}{2 \pi i x}\oint_{L} \frac{\Gamma(qs
)~x^{p ls}}{ (2 \pi)^{\frac{1-p}{2}} p^{pls
-\frac{1}{2}}\Gamma(ls) \Gamma(ls+\frac{1}{p}) \cdots
\Gamma(ls+\frac{p-1}{p})}ds\\
%
%
%
%
&&=\frac{q}{2 \pi i x}\oint_{L} \frac{(2 \pi)^{\frac{1-q}{2}}
q^{qs -\frac{1}{2}}\Gamma(s) \Gamma(s+\frac{1}{q}) \cdots
\Gamma(s+\frac{q-1}{q}) ~x^{p l s}}{ (2 \pi)^{\frac{1-p}{2}} p^{pl
s -\frac{1}{2}}\Gamma(l s) \Gamma(l s+\frac{1}{p})
\cdots \Gamma(l s+\frac{p-1}{p})}ds\\
%
%
&&=\frac{1}{2 \pi i }\oint_{L} \frac{\frac{\sqrt{pq}}{x(2
\pi)^{\frac{q-p}{2}}}~ ls\Gamma(s)
\Gamma\left(s+\frac{1}{q}\right) \cdots
\Gamma\left(s+\frac{q-1}{q}\right) }{ ls\Gamma\left(l s\right)
\Gamma\left(l s+\frac{1}{p}\right) \cdots \Gamma\left(l
s+\frac{p-1}{p}\right)}\left[\frac{p^{pl }}{x^{p l}q^{q}
}\right]^{
-s} ds\\
&&=\frac{l\sqrt{pq}}{x(2 \pi)^{\frac{q-p}{2}}}~\frac{1}{2 \pi i
}\oint_{L} \frac{  \Gamma(s+1) \Gamma\left(s+\frac{1}{q}\right)
\cdots \Gamma\left(s+\frac{q-1}{q}\right)}{  \Gamma\left(l
s+1\right)\Gamma\left(l s+\frac{1}{p}\right) \cdots \Gamma\left(l
s+\frac{p-1}{p}\right)}\left[\frac{p^{pl }}{x^{p l}q^{q}
}\right]^{
-s} ds\\
&&=\frac{l\sqrt{pq}}{x(2
\pi)^{\frac{q-p}{2}}}H^{q,0}_{p,q}\left[\frac{p^{pl }}{x^{p
l}q^{q}}
 \big|^{\left(1,l\right),\left(\frac{1}{p},l\right),\cdots,\left(
\frac{p-1}{p},l\right)}_{(1,1),\left(\frac{1}{q},1
\right),\cdots,\left(\frac{q-1}{q},1\right) }\right]
\end{eqnarray*}
Since we have the following
\begin{eqnarray*}\label{}
&&H^{q,0}_{p,q}\left[\frac{p^{pl }}{x^{pl }q^{q}}
 \big|^{\left(1,l\right),\left(\frac{1}{p},l\right),\cdots,\left(
\frac{p-1}{p},l\right)}_{(1,1),\left(\frac{1}{q},1
\right),\cdots,\left(\frac{q-1}{q},1\right) }\right]\\
&&=\sum^{\infty}_{k_{0}=1}\frac{(-1)^{k_{0}}\Gamma(-k_{0}-\frac{q-1}{q})
\cdots \Gamma(-k_{0}-\frac{2}{q}) \Gamma(-k_{0}-\frac{1}{q}) } {
k_{0}!~\Gamma(- lk_{0}+\frac{pq-plq}{pq})
\Gamma(-lk_{0}+\frac{q-plq}{pq})
\cdots \Gamma(- lk_{0}+\frac{(p-1)q-plq}{pq}) }\left(
\frac{p^{pl}}{x^{pl}q^{q}}\right)^{k_{0}+1}
\end{eqnarray*}
\begin{eqnarray*}
&&+\sum^{\infty}_{k_{1}=0}\frac{(-1)^{k_{1}}\Gamma(-k_{1}+\frac{q-1}{q})
\Gamma(-k_{1}+\frac{1}{q})\cdots \Gamma(-k_{1}+\frac{q-2}{q}) } {
k_{1}!~ \Gamma(-
lk_{1}+\frac{pq-pl}{pq})\Gamma(-lk_{1}+\frac{q-pl}{pq})\cdots
\Gamma(- lk_{1}+\frac{(p-1)q-pl}{pq}) }\left(
\frac{p^{pl}}{x^{pl}q^{q}}\right)^{k_{1}+\frac{1}{q}}\\
&&+\sum^{\infty}_{k_{2}=0}\frac{(-1)^{k_{2}}\Gamma(-k_{2}+\frac{q-2}{q})
\Gamma(-k_{2}-\frac{1}{q})\Gamma(-k_{2}+\frac{1}{q})\cdots
\Gamma(-k_{2}+\frac{q-3}{q}) } { k_{2}!~\Gamma(-
lk_{2}+\frac{pq-2pl}{pq}) \Gamma(-lk_{2}+\frac{q-2pl}{pq})\cdots
\Gamma(- lk_{2}+\frac{(p-1)q-2pl}{pq}) }\left(
\frac{p^{pl}}{x^{pl}q^{q}}\right)^{k_{2}+\frac{2}{q}}\\
&&+\cdots\\
%
&&+\sum^{\infty}_{k_{q-1}=0}\frac{(-1)^{k_{q-1}}\Gamma(-k_{q-1}+\frac{1}{q})
\Gamma(-k_{q-1}-\frac{q-2}{q})\cdots \Gamma(-k_{q-1}-\frac{1}{q})
} { k_{q-1}!~ \Gamma(-
lk_{q-1}+\frac{pq-(q-1)pl}{pq})\Gamma(-lk_{q-1}+\frac{q-(q-1)pl}{pq})\cdots
\Gamma(- lk_{q-1}+\frac{(p-1)q-(q-1)pl}{pq}) }\\&&\times \left(
\frac{p^{pl}}{x^{pl}q^{q}}\right)^{k_{q-1}+\frac{q-1}{q}}
=\left(\frac{p^{lp}}{x^{lp}
q^{q}}\right)~_{q-1}\Psi_{p}\left[-\frac{p^{pl}}{x^{pl}q^{q}}
{\Huge |}^{(-\frac{q-1}{q},-1)\cdots
 (-\frac{2}{q},-1) (-\frac{1}{q},-1)}_{(1-l,-l)
(\frac{q-plq}{pq},-l)\cdots (\frac{(p-1)q-plq}{pq},-l)}\right]\\
%
&&+\left(\frac{p^{l\frac{p}{q}}}{x^{l\frac{p}{q}} q}\right)
   ~_{q-1}\Psi_{p}\left[-\frac{p^{pl}}{x^{pl}q^{q}} {\Huge
|}^{(\frac{q-1}{q},-1)(\frac{1}{q},-1)\cdots (\frac{q-2}{q},-1)
}_{
(\frac{pq-pl}{p},-l)(\frac{q-pl}{pq},-l)\cdots (\frac{(p-1)q-pl}{p},-l)}\right]\\
&&+\left(\frac{p^{2l\frac{p}{q}}}{x^{2l\frac{p}{q}} q^{2}}\right)
   ~_{q-1}\Psi_{p}\left[-\frac{p^{pl}}{x^{pl}q^{q}} {\Huge
|}^{(-\frac{1}{q},-1)(\frac{1}{q},-1)\cdots (\frac{q-3}{q},-1)
(\frac{q-2}{q},-1)}_{
(\frac{q-2pl}{pq},-l)\cdots (\frac{(p-1)q-2pl}{p},-l)(\frac{pq-2pl}{p},-l)}\right]\\
&&+\cdots\\
&&+\left(\frac{p^{(q-1)l\frac{p}{q}}}{x^{(q-1)l\frac{p}{q}}
q^{(q-1)}}\right)
   ~_{q-1}\Psi_{p}\left[-\frac{p^{pl}}{x^{pl}q^{q}} {\Huge
|}^{(\frac{1}{q},-1)(-\frac{(q-2)}{q},-1)(\frac{(q-3)}{q},-1)\cdots
(-\frac{1}{q},-1) }_{
(\frac{pq-(q-1)pl}{p},-l)(\frac{q-(q-1)pl}{pq},-l)\cdots
(\frac{(p-1)q-(q-1)pl}{p},-l)}\right],
\end{eqnarray*}
we get
\begin{eqnarray}
\nonumber&&f_{(\alpha,l)}(x)=\frac{l\sqrt{pq}}{x(2
\pi)^{\frac{q-p}{2}}}~H^{q,0}_{p,q}\left[\frac{p^{pl }}{x^{pl
}q^{q}}
 \big|^{\left(1,l\right),\left(\frac{1}{p},l\right),\cdots,\left(
\frac{p-1}{p},l\right)}_{(1,1),\left(\frac{1}{q},1
\right),\cdots,\left(\frac{q-1}{q},1\right) }\right]\\
%
%
\nonumber &&=\frac{l\sqrt{pq}}{x(2
\pi)^{\frac{q-p}{2}}}\left(\frac{p^{lp}}{x^{lp}
q^{q}}\right)~_{q-1}\Psi_{p}\left[-\frac{p^{pl}}{x^{pl}q^{q}}
{\Huge |}^{(-\frac{q-1}{q},-1)\cdots
 (-\frac{2}{q},-1) (-\frac{1}{q},-1) }_{(1-l,-l)
(\frac{q-plq}{pq},-l)\cdots (\frac{(p-1)q-plq}{pq},-l)}\right]\\
\nonumber &&+\frac{l\sqrt{pq}}{x(2
\pi)^{\frac{q-p}{2}}}\sum^{q-1}_{j=1}
\left(\frac{p^{jl\frac{p}{q}}}{x^{jl\frac{p}{q}} q^{j}}\right)
~_{q-1}\Psi_{p}\left[-\frac{p^{pl}}{x^{pl}q^{q}} {\Huge
|}^{\left(1-\left[\frac{j}{q}\right]_{q},-1\right)\left(\frac{2}{q}-\left[\frac{j}{q}\right]_{2},-1\right)
\left(\frac{3}{q}-\left[\frac{j}{q}\right]_{3},-1\right)\cdots
\left(\frac{q-1}{q}-\left[\frac{j}{q}\right]_{q-1},-1\right)
 }_{\left(\frac{pq-j
pl}{pq},-l\right)\left(\frac{q-j pl}{pq},-l\right)\cdots
\left(\frac{(p-1)q-j pl}{pq},-l\right)}\right]\label{}
\end{eqnarray}
where $\left[\frac{j}{q}\right]_{n}$
is  the L\'{e}vy jump function.\\ Note that the series are
absolutely convergent satisfying the condition
$\sum_{j=1}^{q}B_{j}-\sum_{i=1}^{p}A_{i}>-1 \Rightarrow
-lp+q-1>-1\Rightarrow q>lp\Rightarrow  \frac{q}{p}>l \Rightarrow
\frac{q}{p}> \frac{p}{q}\left(\frac{p}{q}
\right)^{\frac{l_{2}}{l_{1}}-1} \Rightarrow \frac{q}{p}>1>
\left(\frac{p}{q}
\right)^{\frac{l_{2}}{l_{1}}}=\alpha$ since $\alpha<1$.\\
If $l$ belongs to the set $\{2,3,4,\cdots\}$, then we have
\begin{eqnarray}
&&f_{(\alpha,l)}(x)=\frac{\sqrt{pq}}{x(2
\pi)^{\frac{q-p}{2}}}H^{q,0}_{p,q}\left[\frac{p^{pl }}{x^{p
l}q^{q}}
 \big|^{\left(0,l\right),\left(\frac{1}{p},l\right),\cdots,\left(
\frac{p-1}{p},l\right)}_{(0,1),\left(\frac{1}{q},1
\right),\cdots,\left(\frac{q-1}{q},1\right) }\right]\label{h5}\\
\nonumber &&=\frac{\sqrt{pq}}{x(2 \pi)^{\frac{q-p}{2}}}\frac{1}{2
\pi i }\oint_{L} \frac{~ \Gamma(s)
\Gamma\left(s+\frac{1}{q}\right) \cdots
\Gamma\left(s+\frac{q-1}{q}\right) }{ \Gamma\left(l s\right)
\Gamma\left(l s+\frac{1}{p}\right) \cdots \Gamma\left(l
s+\frac{p-1}{p}\right)}\left[\frac{p^{pl }}{x^{p l}q^{q}
}\right]^{ -s} ds
\end{eqnarray}
We use
$$\Gamma(l s)=(2 \pi)^{\frac{1-l}{2}} l^{ls -\frac{1}{2}}\Gamma(s)
\Gamma\left(s+\frac{1}{l}\right) \cdots
\Gamma\left(s+\frac{l-1}{l}\right).$$ Hence we have
\begin{eqnarray*}
%
&&=\frac{\sqrt{pq}}{x(2 \pi)^{\frac{q-p}{2}}}\frac{1}{2 \pi i
}\oint_{L} \frac{~ \Gamma(s) \Gamma\left(s+\frac{1}{q}\right)
\cdots \Gamma\left(s+\frac{q-1}{q}\right) }{(2
\pi)^{\frac{1-l}{2}} l^{ls -\frac{1}{2}}\Gamma(s)
\Gamma\left(s+\frac{1}{l}\right) \cdots
\Gamma\left(s+\frac{l-1}{l}\right) \Gamma\left(l
s+\frac{1}{p}\right) \cdots \Gamma\left(l
s+\frac{p-1}{p}\right)}\left[\frac{p^{pl }}{x^{p l}q^{q}
}\right]^{ -s} ds\\
&&=\frac{\sqrt{pql}}{x(2 \pi)^{\frac{q+1-p-l}{2}}}\frac{1}{2 \pi i
}\oint_{L} \frac{ \Gamma\left(s+\frac{1}{q}\right) \cdots
\Gamma\left(s+\frac{q-1}{q}\right) }{
\Gamma\left(s+\frac{1}{l}\right) \cdots
\Gamma\left(s+\frac{l-1}{l}\right) \Gamma\left(l
s+\frac{1}{p}\right) \cdots \Gamma\left(l
s+\frac{p-1}{p}\right)}\left[\frac{p^{pl }l^{l}}{x^{p l}q^{q}
}\right]^{ -s} ds\\
&&=\frac{\sqrt{pql}}{x(2
\pi)^{\frac{q+1-p-l}{2}}}H^{q-1,0}_{p+l-2,q-1}\left[\frac{p^{pl
}l^{l}}{x^{p l}q^{q}}
 \big|^{\left(\frac{1}{l},1\right),\cdots,\left(\frac{l-1}{l},1\right),\left(\frac{1}{p},l\right),\cdots,\left(
\frac{p-1}{p},l\right)}_{\left(\frac{1}{q},1
\right),\cdots,\left(\frac{q-1}{q},1\right) }\right]
\end{eqnarray*}
Note that for $l\in\{2,3,4,\cdots\}$, $q$ should be of the form
$kp$, $k=4,5,6,\cdots$ and $l_{2}<l_{1}$.
$\sum_{j=1}^{q}B_{j}-\sum_{i=1}^{p}A_{i}>0 $ means $q-1 -(l-1)
-l(p-1)  >0 \Rightarrow q>lp \Rightarrow kp>lp \Rightarrow k>l$.
But this condition is always satisfied since
$l=(k)^{1-\frac{l_{2}}{l_{1}}}$.\\

If we put $l=1$ in \eqref{h5}, then
\begin{eqnarray}
\nonumber&&f_{(\alpha,1)}(x)=\frac{\sqrt{pq}}{x(2
\pi)^{\frac{q-p}{2}}}~H^{q,0}_{p,q}\left[\frac{p^{p }}{x^{p
}q^{q}}
 \big|^{\left(0,1\right),\left(\frac{1}{p},1\right),\cdots,\left(
\frac{p-1}{p},1\right)}_{(0,1),\left(\frac{1}{q},1
\right),\cdots,\left(\frac{q-1}{q},1\right) }\right]
=\frac{\sqrt{pq}}{x(2
\pi)^{\frac{q-p}{2}}}~H^{q-1,0}_{p-1,q-1}\left[\frac{p^{p }}{x^{p
}q^{q}}
 \big|^{\left(\frac{1}{p},1\right),\cdots,\left(
\frac{p-1}{p},1\right)}_{\left(\frac{1}{q},1
\right),\cdots,\left(\frac{q-1}{q},1\right) }\right]\\
%
\nonumber&&=\frac{\sqrt{pq}}{x(2
\pi)^{\frac{q-p}{2}}}\sum^{q-1}_{j=1}
\left(\frac{p^{j\frac{p}{q}}}{x^{j\frac{p}{q}} q^{j}}\right)
~_{q-2}\Psi_{p-1}\left[-\frac{p^{p}}{x^{p}q^{q}} {\Huge
|}^{\left(\frac{2}{q}-\left[\frac{j}{q}\right]_{2},-1\right)
\left(\frac{3}{q}-\left[\frac{j}{q}\right]_{3},-1\right)\cdots
\left(\frac{q-1}{q}-\left[\frac{j}{q}\right]_{q-1},-1\right) }_{
\left(\frac{q-j p}{pq},-1\right)\cdots \left(\frac{(p-1)q-j
p}{pq},-1\right)}\right]\\
 \nonumber&& \mbox{by applying the formula $\Gamma(\beta+1-v)=
\frac{(-1)^{v}\Gamma(\beta+1)}{(1-(\beta+1))_{v}}$ in Lemma 1.1}\\
 \nonumber&&=\frac{\sqrt{pq}}{x(2
\pi)^{\frac{q-p}{2}}}\sum^{q-1}_{j=1}
\left(\frac{p^{j\frac{p}{q}}}{x^{j\frac{p}{q}} q^{j}}\right)
\frac{
\prod_{i_{1}=2}^{q-1}\Gamma\left(\frac{i_{1}}{q}-\left[\frac{j}{q}\right]_{i_{1}}\right)
} {\prod_{i_{2}=1}^{p-1}\Gamma\left(\frac{i_{2}q-j p}{pq}\right)}
\\
\nonumber&&\times
~_{p-1}F_{q-2}\left[(-1)^{q-p}\frac{p^{p}}{x^{p}q^{q}} {\Huge
|}_{\left(1-\frac{2}{q}+\left[\frac{j}{q}\right]_{2}\right)
\left(1-\frac{3}{q}+\left[\frac{j}{q}\right]_{3}\right)\cdots
\left(1-\frac{q-1}{q}+\left[\frac{j}{q}\right]_{q-1}\right) }^{
\left(1-\frac{q-j p}{pq}\right)\cdots \left(1-\frac{(p-1)q-j
p}{pq}\right)}\right]
\end{eqnarray}
 For the case of $l=1$, the condition for their convergence is
trivial since $p<q$.\\


Some special cases will be given.
We will use \eqref{h2} to show some known results.\\
\begin{ex}
When $p=1,~q=4,~ \alpha=\frac{1}{4}$, we have
\begin{eqnarray*}\label{}
&&f_{(\frac{1}{4},1)}(x)=\frac{2}{x(2
\pi)^{\frac{3}{2}}}~H^{3,0}_{0,3}\left[\frac{1}{4^{4}x}
 \big|^{-}_{\left(\frac{1}{4},1
\right),\left(\frac{2}{4},1 \right),\left(\frac{3}{4},1\right)
}\right]\\
&&=\frac{2}{  x(2 \pi)^{\frac{3}{2}}} \sum^{3}_{j=1}
\left(\frac{1}{x^{j\frac{1}{4}} 4^{j}}\right) \frac{
\prod_{i_{1}=2}^{3}\Gamma\left(\frac{i_{1}}{4}-\left[\frac{j}{4}\right]_{i_{1}}\right)
} {\prod_{i_{2}=1}^{0}\Gamma\left(\frac{i_{2}4-j }{4}\right)}
~_{0}F_{2}\left[(-1)^{3}\frac{1}{4^{4}x} {\Huge
|}_{\left(1-\frac{2}{4}+\left[\frac{j}{4}\right]_{2}\right)
\left(1-\frac{3}{4}+\left[\frac{j}{4}\right]_{3}\right) }^{
-}\right]\\
&&=\frac{1}{2x^{\frac{5}{4}}(2 \pi)^{\frac{3}{2}}}
\Gamma\left(\frac{1}{4}\right) \Gamma\left(\frac{1}{2}\right)
~_{0}F_{2} \left[\frac{-1}{256x}
\big|^{-}_{\left(\frac{3}{4}\right)
\left(\frac{1}{2}\right)}\right]
+\frac{1}{8x^{\frac{3}{2}}(2 \pi)^{\frac{3}{2}}}
\Gamma\left(-\frac{1}{4}\right) \Gamma\left(\frac{1}{4}\right)
~_{0}F_{2} \left[\frac{-1}{256x}
\big|^{-}_{\left(\frac{3}{4}\right)
\left(\frac{5}{4}\right)}\right]
\end{eqnarray*}
\begin{eqnarray*}
&&+\frac{1}{32x^{\frac{7}{4}}(2 \pi)^{\frac{3}{2}}}
\Gamma\left(-\frac{1}{4}\right) \Gamma\left(-\frac{1}{2}\right)
~_{0}F_{2} \left[\frac{-1}{256x}
\big|^{-}_{\left(\frac{5}{4}\right)
\left(\frac{3}{2}\right)}\right]\\
&&=\frac{2}{(2 \pi)^{\frac{3}{2}}x}\{ \frac{2}{x^{\frac{1}{4}}}
\Gamma\left(1+\frac{1}{4}\right) \Gamma\left(1+\frac{2}{4}\right)
~_{0}F_{2} \left[\frac{-1}{256x}
\big|^{-}_{\left(\frac{3}{4}\right)
\left(\frac{2}{4}\right)}\right]\\
&&-\frac{1}{x^{\frac{1}{2}}} \Gamma\left(\frac{3}{4}\right)
\Gamma\left(1+\frac{1}{4}\right) ~_{0}F_{2} \left[\frac{-1}{256x}
\big|^{-}_{\left(\frac{3}{4}\right)
\left(1+\frac{1}{4}\right)}\right]
+\frac{1}{8x^{\frac{3}{4}}} \Gamma\left(\frac{2}{4}\right)
\Gamma\left(\frac{3}{4}\right) ~_{0}F_{2} \left[\frac{-1}{256x}
\big|^{-}_{\left(1+\frac{1}{4}\right)
\left(1+\frac{2}{4}\right)}\right]\}.
\end{eqnarray*}
\end{ex}

We will use \eqref{h1} to show some new results.\\
\begin{ex}
Let  $p=1,~q=2,~l_{1}=2,~l_{2}=1, ~\alpha=\frac{1}{\sqrt{2}},
~l=\sqrt{2}$, then
\begin{eqnarray*}
\nonumber&&f_{(\frac{1}{\sqrt{2}},\sqrt{2})}(x)=\frac{\sqrt{2}}{x\sqrt{\pi}}
~H^{2,0}_{1,2} \left[\frac{1}{4x^{\sqrt{2} }}
\big|^{\left(1,\sqrt{2}\right)}_{(1,1),\left(\frac{1}{2},1
\right) }\right]
\\\nonumber&&=\frac{\sqrt{2}}{x\sqrt{\pi}} \left( \frac{1}{4 x^{\sqrt{2}}}
~_{1}\Psi_{1}\left[-\frac{1}{4x^{\sqrt{2}}} {\Huge
|}^{\left(-\frac{1}{2},-1\right)
 }_{\left(1-\sqrt{2},-\sqrt{2}\right) }\right]+\frac{1}{2
x^{\frac{1}{\sqrt{2}}}}
~_{1}\Psi_{1}\left[-\frac{1}{4x^{\sqrt{2}}} {\Huge
|}^{\left(\frac{1}{2},-1\right)
 }_{\left(1-\frac{1}{\sqrt{2}},-\sqrt{2}\right) }\right]\right). 
%
\end{eqnarray*}
\end{ex}

\section{A process to generate more representations}

For the same rational $\alpha$, more representations of one
distribution can be generated by using \eqref{h2} and \eqref{h3}.
\begin{ex}
When $p=1,~q=2,~ \alpha=\frac{1}{2}$ in \eqref{h2}, we have
\begin{eqnarray*}\label{}
&&f_{(\frac{1}{2},1)}(x)=\frac{1}{x \sqrt{\pi}}
~H^{1,0}_{0,1}\left[\frac{1}{4x}
 \big|^{-}_{\left(\frac{1}{2},1
\right) }\right]=\frac{1}{ 2 x^{\frac{3}{2}} \sqrt{\pi}   }
~_{0}F_{0}\left[\frac{-1}{4x} {\Huge |}_{- }^{ -}\right]=\frac{1}{
2\sqrt{\pi} x^{\frac{3}{2}}    }exp\left(\frac{-1}{4x} \right)
\end{eqnarray*}
\end{ex}

\begin{ex}
Let  $p=1,~q=4,~l_{1}=2,~l_{2}=1,
~\alpha=\frac{1}{4}^{\frac{1}{2}}=\frac{1}{2}, ~l=2$ in
\eqref{h3}, then
we have
\begin{eqnarray*}\label{}
&&f_{(\frac{1}{2},2)}(x)=\frac{\sqrt{8}}{2 \pi
x}H^{3,0}_{1,3}\left[\frac{4}{x^{2}4^{4}}
  \big|^{\left(\frac{1}{2},1\right)}_{\left(\frac{1}{4},1
\right),\left(\frac{2}{4},1\right),\left(\frac{3}{4},1\right)
}\right]=\frac{1}{2 \pi x^{\frac{3}{2}}}\sum_{v=0}^{\infty}
\frac{\Gamma(-v+\frac{1}{2})}{v!}
\left(\frac{-1}{4^{3}x^{2}}\right)^{v}\\
&&+\frac{1}{2^{4} \pi x^{\frac{5}{2}}}\sum_{v=0}^{\infty}
\frac{\Gamma(-v-\frac{1}{2})}{v!}
\left(\frac{-1}{4^{3}x^{2}}\right)^{v}=\frac{1}{2\sqrt{\pi}
x^{\frac{3}{2}}}\left(
~_{0}F_{1}\left[\frac{1}{4^{3}x^{2}}|^{-}_{\left(\frac{1}{2}\right)}\right]-
\frac{1}{4x}
~_{0}F_{1}\left[\frac{1}{4^{3}x^{2}}|^{-}_{\left(\frac{3}{2}\right)}
\right]\right)
\end{eqnarray*}
\end{ex} Note that $exp\left(\frac{-1}{4x} \right) =
~_{0}F_{1}\left[\frac{1}{4^{3}x^{2}}|^{-}_{\left(\frac{1}{2}\right)}\right]-
\frac{1}{4x}
~_{0}F_{1}\left[\frac{1}{4^{3}x^{2}}|^{-}_{\left(\frac{3}{2}\right)}
\right]$ and by setting
$\alpha=\left(\frac{2}{8}\right)^{\frac{1}{2}}=\frac{1}{2}$ in
\eqref{h3}, another representation can be born in the sum of
faster convergent series.

\section{L\'{e}vy smashing on the family of gamma density functions and L\'{e}vy-smashed gamma stochastic process}

Mellin transform of a density function in statistics shows its
statistical structure and this technique can be used as a tool to
blend two independently distributed random variables. In this
section, we show what L\'{e}vy effect could be and how we should
understand it. To start with, consider the L\'{e}vy density
function $\frac{1}{2 \pi i} \oint_{L}
\frac{\Gamma(\frac{1}{\alpha}-\frac{s}{\alpha})}{\alpha
\Gamma(1-s)}x^{-s}ds$ and the gamma density functions
$\frac{x^{\gamma-1}}{\Gamma(\gamma)}e^{-x}$ where
$0<\alpha<1,~0<\gamma$ and $0<x<\infty$. $\frac{1}{2 \pi i}
\oint_{L} \frac{\Gamma(\frac{1}{\alpha}-\frac{s}{\alpha})}{\alpha
\Gamma(1-s)}x^{-s}ds$ is the one-sided L\'{e}vy density function
 found in \cite{ma1} constructed in a different way in \cite{h1}.
 We will use the Mellin transformation of the form
 $\int_{t}h_{1}(\frac{x}{t})h_{2}(t)\frac{dt}{t}$, where $h_{1}(x)$ and
 $h_{2}(x)$ are certain density functions. Then we have
\begin{eqnarray}\label{}
\nonumber &&f_{(\alpha)}(x)=\int^{\infty}_{0} \frac{1}{2 \pi i}
\oint_{L} \frac{\Gamma(\frac{1}{\alpha}-\frac{s}{\alpha})}{\alpha
\Gamma(1-s)}x^{-s}t^{s}ds
~\frac{t^{\gamma-1}}{\Gamma(\gamma)}e^{-t}\frac{dt}{t}
=\frac{1}{2 \pi i} \oint_{L}
\frac{\Gamma(\frac{1}{\alpha}-\frac{s}{\alpha})}{\alpha
\Gamma(1-s)}x^{-s}\int^{\infty}_{0}
\frac{t^{s+\gamma-1-1}}{\Gamma(\gamma)}e^{-t}dt ds\\
\nonumber&&=\frac{1}{2 \pi i} \oint_{L}
\frac{\Gamma(\frac{1}{\alpha}-\frac{s}{\alpha})}{\alpha
\Gamma(1-s)} \frac{\Gamma(s+\gamma-1)}{\Gamma(\gamma)} x^{-s} ds
~~\mbox{by transformation } s=1-\alpha s_{1},\\
&&=\frac{1}{2 \pi i} \oint_{L} \frac{\Gamma(s_{1})}{ \Gamma(\alpha
s_{1})}\frac{\Gamma(\gamma-\alpha s_{1})}{\Gamma(\gamma)}x^{\alpha
s_{1}-1} ds_{1} \label{h6}
\end{eqnarray}
and its Laplace transform
\begin{eqnarray}\label{}
\nonumber &&L_{f}(y)=\int_{0}^{\infty}e^{-yx} \frac{1}{2 \pi i}
\oint_{L} \frac{\Gamma(s)}{ \Gamma(\alpha
s)}\frac{\Gamma(\gamma-\alpha s)}{\Gamma(\gamma)} x^{\alpha s-1}
ds dx
=\frac{1}{2 \pi i} \oint_{L} \frac{\Gamma(s)\Gamma(\gamma-\alpha
s)}{\Gamma(\gamma)} y^{-\alpha s} ds\\
&&=\sum_{k=0}^{\infty} \frac{(-1)^{k} \Gamma(\gamma+\alpha
k)}{k!~\Gamma(\gamma)} \left(y\right)^{\alpha k}.
\end{eqnarray}
 When
$\alpha=1$ in \eqref{h6}, \eqref{h6} becomes gamma density
$\frac{x^{\gamma-1}}{\Gamma(\gamma)}e^{-x}$. So $f_{(1)}(x)$ is a
one-sided function concentrated at $x=1$ for $\mathbb{R}^{+}$.
To understand its effect, put $\alpha=\frac{1}{2}$ in \eqref{h6},
then we have
\begin{eqnarray}\label{}
&&\frac{1}{2 \pi i} \oint_{L} \frac{\Gamma(s_{1})}{
\Gamma(\frac{1}{2} s_{1})}\frac{\Gamma(\gamma-\frac{1}{2}
s_{1})}{\Gamma(\gamma)} x^{\frac{1}{2} s_{1}-1} ds_{1}, s_{1}=2s\\
&&=\frac{1}{2 \pi i} \oint_{L} \frac{\Gamma(2s)}{ \Gamma(
s)}\frac{\Gamma(\gamma-s)}{\Gamma(\gamma)} x^{ s-1} 2ds
=\frac{1}{2 \pi i} \oint_{L} \frac{2^{2s-1}\Gamma(s)
\Gamma(s+\frac{1}{2})}{\sqrt{\pi} \Gamma(
s)}\frac{\Gamma(\gamma-s)}{\Gamma(\gamma)} x^{ s-1} 2ds\\
&&=\sum_{k=0}^{\infty} \frac{(-1)^{k}
\Gamma(\gamma+\frac{1}{2}+k)} {k! \sqrt{\pi}   \Gamma(\gamma)}
x^{-k-\frac{3}{2}}4^{-k-\frac{1}{2}}
=\frac{\Gamma(\gamma+\frac{1}{2})}{2\sqrt{\pi}\Gamma(\gamma)x^{\frac{3}{2}}}
\left(1+\frac{1}{4x}\right)^{-\gamma-\frac{1}{2}}\\
&&=\frac{4\Gamma(\gamma+\frac{1}{2})(4x)^{\gamma-1}}{\sqrt{\pi}\Gamma(\gamma)}
\left(1+4x\right)^{-\gamma-\frac{1}{2}}
\end{eqnarray}
Fig. 1 shows the impact on the family of gamma functions.\\
\begin{tabular}{|c|c|c|c|}
  \hline
  &gamma family&&L\'{e}vy smashed gamma family\\
  \hline
  (a) $~\gamma=1$ & $e^{-x}$ &$\leftrightarrow $ &$2(1+4x)^{-3/2}$  \\
  \hline
  (b) $~\gamma=2$& $xe^{-x}$ & $\leftrightarrow$ & $12x(1+4x)^{-5/2}$  \\
\hline
  (c) $~\gamma=3$& $\frac{1}{2}x^{2}e^{-x}$ &$\leftrightarrow $ & $60x^{2}(1+4x)^{-7/2}$  \\
\hline
  (d) $~\gamma=4$ & $\frac{1}{3!}x^{3}e^{-x}$ &$\leftrightarrow$  & $280x^{3}(1+4x)^{-9/2}$  \\
\hline
  \end{tabular}
\begin{center}
\resizebox{7.1cm}{6cm}{\includegraphics{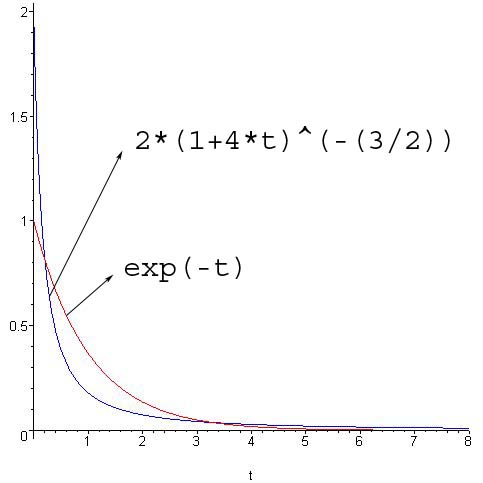}}(a)
\resizebox{7.1cm}{6cm}{\includegraphics{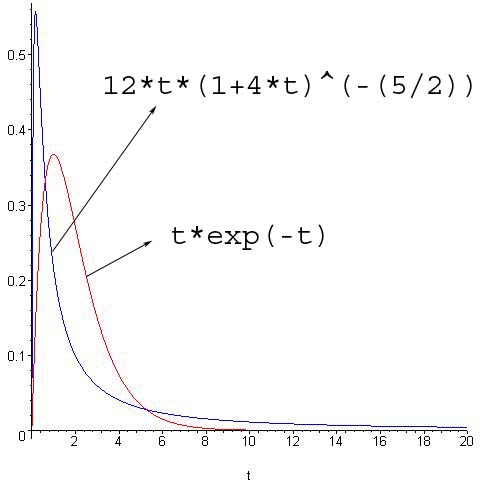}}(b)\\
\resizebox{7cm}{6cm}{\includegraphics{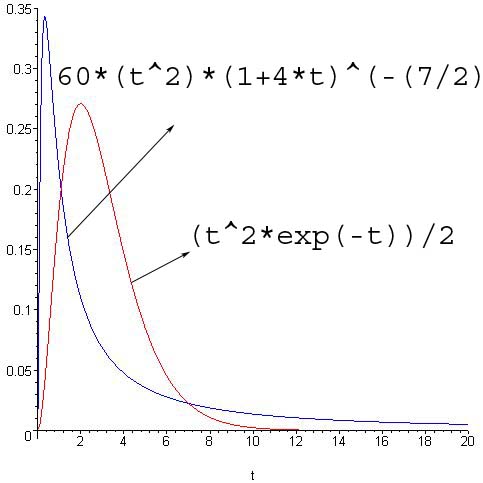}} (c)
\resizebox{7cm}{6cm}{\includegraphics{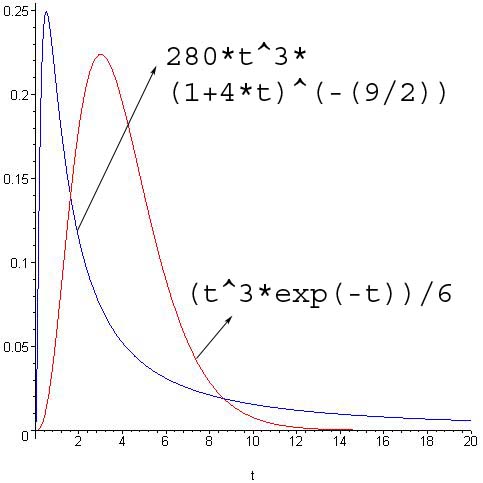}} (d)\\
Figure 1.
\end{center}

$f_{(\alpha)}(x)=\frac{1}{2 \pi i} \oint_{L} \frac{\Gamma(s)}{
\Gamma(\alpha s)}\frac{\Gamma(\gamma-\alpha
s)}{\Gamma(\gamma)}x^{\alpha s-1} ds$ is absolutely convergent
series for all $x$ since $\mu=\alpha-\alpha+1=1>0$, see
\cite{mhs2}. From the fig. 1, $f_{(\alpha)}(x)$ will be called as
L\'{e}vy-smashed gamma density functions especially when the
parameter $0<\alpha<1$. Note that $\alpha$ can be any positive real number.\\

 The stochastic process {$X(t), t> 0$} with
$X(0)= 0$ and having stationary and independent increments, where
$X(t)$ has the density function $\frac{1}{2 \pi i} \oint_{L}
\frac{\Gamma(s)}{ \Gamma(\alpha s)}\frac{\Gamma(t-\alpha
s)}{\Gamma(t)} x^{\alpha s-1} ds,~ 0 < \alpha\leq 1,$ will be
called L\'{e}vy-smashed gamma stochastic process. The
L\'{e}vy-smashed gamma stochastic process $X(t)$ has the
distribution function,for $t > 0,~~ 0 < \alpha < 1$,
$F_{(\alpha,t)}(x)=\frac{1}{2 \pi i} \oint_{L} \frac{\Gamma(s)}{
\Gamma(1+\alpha s)}\frac{\Gamma(t-\alpha s)}{\Gamma(t)}x^{\alpha
s} ds
$. From \cite{f1}
and \cite{p1}, we can prove that the L\'{e}vy-smashed gamma
distribution with parameter $\alpha$ is attracted to the stable
distribution with exponent $\alpha,~~ 0 < \alpha < 1$. Namely,
\begin{eqnarray}\label{}
\nonumber&&\lim_{n\rightarrow \infty}
L_{f}(\frac{y}{n})=\lim_{n\rightarrow
\infty}\int_{0}^{\infty}e^{-yx} \frac{1}{2 \pi i} \oint_{L}
\frac{\Gamma(s)}{ \Gamma(\alpha s)}\frac{\Gamma(n-\alpha
s)}{\Gamma(n)}n^{\alpha s} x^{\alpha s-1} ds dx\\
\nonumber&&=\lim_{n\rightarrow \infty}\frac{1}{2 \pi i} \oint_{L}
\frac{\Gamma(s)\Gamma(n-\alpha s)}{\Gamma(n)}n^{\alpha s}
y^{-\alpha s} ds
=\lim_{n\rightarrow \infty}\sum_{k=0}^{\infty} \frac{(-1)^{k}
\Gamma(n+\alpha k)}{k!~ \Gamma(n)}
\left(\frac{y}{n}\right)^{\alpha k}=e^{-y^{\alpha}}
\end{eqnarray}

\section{Remarks}

In \cite{mpg1}, they consider the signalling problem for the
standard diffusion equation $\frac{\partial}{\partial t} u(x,t) =D
\frac{\partial^{2}}{\partial x^{2}} u(x,t)$ with the conditions
$u(x,0^{+})=0~x>0; u(0^{+},t)=h(t), u(+\infty,t)=0,~ t>0.$  And
they say "$\cdots$ Then the solution is as follows
$u(x,t)=\int_{0}^{t}\mathcal{G}_{s}^{d}(x,\tau)h(t-\tau)d\tau$,
where $\mathcal{G}_{s}^{d}(x,t)=\frac{x}{2\sqrt{\pi
D}}t^{-\frac{3}{2}}exp{-\frac{x^{2}}{4Dt}}$. Here
$\mathcal{G}_{s}^{d}(x,t)$ represents the fundamental solution (or
Green function) of the signalling problem, since it corresponds to
$h(t)=\delta(t)$. We note that
\begin{equation}\label{h8}
\mathcal{G}_{s}^{d}(x,t)=p_{LS}(t;\mu):=\frac{\sqrt{\mu}}{\sqrt{2\pi}t^{\frac{3}{2}}}
exp(-\frac{\mu}{2t}),~t\geq0,~\mu=\frac{x^{2}}{2D} \end{equation}
where $p_{LS}(t;\mu)$ denotes the one-sided L\'{e}vy-Smirnov pdf
spread out over all non-negative $t$ (the time variable). The
L\'{e}vy-Smirnov pdf has all moments of integer order infinite,
since it decays at infinity as $t^{-\frac{3}{2}}$. However, we
note that the absolute moments of real order $\nu$ are finite only
if $0\leq \nu <\frac{1}{2}$. In particular, for this pdf the mean
is infinite, for which we can take the median as expected value.
From $\mathcal{P}_{LS}(t_{med};\mu)=\frac{1}{2}$, it turns out
that $t_{med}\approx 2\mu$, since the complementary error function
gets the value $\frac{1}{2}$ as its argument is approximatively
$\frac{1}{2}$. $\cdots$".\\

With the inspiration from the above paragraph, take the L\'{e}vy
density function with the parameter $\alpha=\frac{1}{2}$, then the
density function is well known to be $\frac{1}{2\sqrt{\pi}
t^{\frac{3}{2}}} exp(-\frac{1}{4t})$. Put this in the mellin
convolution formula $\int_{t}h_{1}(\frac{x}{t})
h_{2}(t)\frac{dt}{t}$, then it becomes $\int_{0}^{\infty}
\frac{\sqrt{t}}{2\sqrt{\pi}x^{\frac{3}{2}}}
exp\left(-\frac{t}{4x}\right)~h_{2}(t)dt$.
$\frac{\sqrt{t}}{2\sqrt{\pi}x^{\frac{3}{2}}}
exp\left(-\frac{t}{4x}\right)$ has the same form with \eqref{h8}
when $t=\mu$. Therefore the cases of L\'{e}vy smashing treated in
section 4 can be thought of as superstatistics in Physics and
Bayesian statistical analysis, subordination in statistics,
namely,
\begin{eqnarray*}\label{}
&&f_{(\frac{1}{2})}(x)=\int^{\infty}_{0} \frac{1}{2 \pi i}
\oint_{L} \frac{\Gamma(2-2s)}{\frac{1}{2}
\Gamma(1-s)}x^{-s}y^{s}ds
~\frac{y^{t-1}}{\Gamma(t)}e^{-y}\frac{dy}{y}=\int_{0}^{\infty}
\frac{\sqrt{y}}{2\sqrt{\pi}x^{\frac{3}{2}}}
exp\left(-\frac{y}{4x}\right)~\frac{y^{t-1}}{\Gamma(t)}e^{-y}dy\\
&&=\int_{0}^{\infty} \frac{\sqrt{y}}{2\sqrt{\pi}x^{\frac{3}{2}}}
\sum_{k=0}^{\infty} \frac{(-1)^{k}y^{k}}{k!~(4x)^{k}}
 ~\frac{y^{t-1}}{\Gamma(t)}e^{-y}dy
 = \frac{1}{2\sqrt{\pi}\Gamma(t)x^{\frac{3}{2}}}
\sum_{k=0}^{\infty} \frac{(-1)^{k}}{k!~(4x)^{k}}
 ~\int_{0}^{\infty}y^{t-1+k+\frac{1}{2}}e^{-y}dy\\
 &&\frac{1}{2\sqrt{\pi}\Gamma(t)x^{\frac{3}{2}}}
\sum_{k=0}^{\infty}
\frac{(-1)^{k}\Gamma(t+k+\frac{1}{2})}{k!~(4x)^{k}}
=\frac{4\Gamma(t+\frac{1}{2})(4x)^{t-1}}{\sqrt{\pi}\Gamma(t)}
\left(1+4x\right)^{-t-\frac{1}{2}}
\end{eqnarray*}
But in this paper, the concept of L\'{e}vy smashing is totally
different from the point of view of superstatistics in Physics and
Bayesian statistical analysis, subordination in statistics.

 \textbf{\large Acknowledgement}\\ \textsf{The author would like to
thank to the Department of Science and Technology, Government of
India, New Delhi, for the financial assistance under Project No.
SR/S4/MS:287/05 and the Centre for Mathematical Sciences for
providing all facilities.}

\end{document}